# Representation Of Level Paths Of An Analytic Function

Kerry M. Soileau

June 29, 2006

ABSTRACT

We find an arc-parameterization of the contour on which an given analytic function has constant modulus. This contour is seen to satisfy a differential equation which we explicitly give.

By a <u>level path</u> we mean a differentiable function (path function) which maps an interval on the real line into the complex plane, with the property that on its image, a given non-constant analytic function has constant modulus. So that the level path is a parameterization according to arc length, we require that its derivative has unit modulus everywhere on its interval of definition.

<u>Proposition 1</u>: If $f$ is analytic on an open set $\Re$, and $z_0 \in \Re$ satisfies $f(z_0) \neq 0$ and $f'(z_0) \neq 0$, and for some interval $A \subseteq \mathbb{R}$ and some set $B \subseteq \{z \in \Re; f(z)f'(z) \neq 0\}$ there exists a differentiable function $p: A \to B$ satisfying

$$p'(s) = \frac{\left| f'(p(s)) \overline{f(p(s))} \right| i}{f'(p(s)) \overline{f(p(s))}}$$

and $p(s_0) = z_0$ for some $s_0 \in A$, then $|f(z)| = |f(z_0)|$ for all $z \in p(A)$. Further, $p$ is a path function which parameterizes $p(A)$ according to arc length such that $p(s_0) = z_0$.

<u>Proof</u>: Suppose there exists $p: A \to B$ satisfying

$$p'(s) = \frac{\left| f'(p(s)) \overline{f(p(s))} \right| i}{f'(p(s)) \overline{f(p(s))}}$$

and $p(s_0) = z_0$ for some $s_0 \in A$. Note that for all $s \in A$,

$$\frac{d}{ds}\left(\left|f(p(s))\right|^2\right)$$
$$=\frac{d}{ds}\left(f(p(s))\overline{f(p(s))}\right)=f'(p(s))p'(s)\overline{f(p(s))}+f(p(s))\overline{f'(p(s))}\overline{p'(s)}=0$$

Next, $f$ is clearly continuous on $p(A)$, hence so is $|f|$, thus $|f(z)|$ is a (positive) constant on $p(A)$. Finally, note that $|p'(s)|=\left|\dfrac{\left|f'(p(s))\overline{f(p(s))}\right|i}{f'(p(s))\overline{f(p(s))}}\right|=1$, so $p$ is a path which parameterizes $p(A)$ according to arc length such that $p(s_0)=z_0$, Q.E.D.[1] Thus there exists a path $p:A\to B$ such that $p(A)\subseteq\Re$ with $z_0\in p(A)$ such that $|f(z)|=|f(z_0)|$ for all $z\in p(A)$.

<u>Proposition 2</u>: If $g(z)$ exists such that it is analytic over $\{z\in\Re; f(z)f'(z)\neq 0\}$ and satisfies $f(z)=\exp g(z)$ there, then the level path of $f(z)$ satisfies

$$p'(s)=i\frac{\overline{g'(p(s))}}{|g'(p(s))|}.$$

<u>Proof</u>: If $f(z)=\exp g(z)$, then $f'(z)=g'(z)f(z)$, hence

$$p'(s)=\frac{\left|f'(p(s))\overline{f(p(s))}\right|i}{f'(p(s))\overline{f(p(s))}}=\frac{\left|g'(p(s))f(p(s))\overline{f(p(s))}\right|i}{g'(p(s))f(p(s))\overline{f(p(s))}}=i\frac{\overline{g'(p(s))}}{|g'(p(s))|},\ \text{Q.E.D.}$$

---

[1] The alternative definition $p'(s)=-\dfrac{\left|f'(p(s))\overline{f(p(s))}\right|i}{f'(p(s))\overline{f(p(s))}}$ would have sufficed, and corresponds to a parameterization of $p(A)$ in the opposite direction.

In the following we will write $c = |f(z_0)|$ for ease of notation, so that $|f(z)| = c$ for all $z \in p(A)$. Using this notation, we may write

$$p'(s) = \frac{|f'(p(s))|ci}{f'(p(s))\overline{f(p(s))}}.$$

After squaring both sides and some rearrangement we get that $f'(p(s))\overline{f(p(s))}^2 p'(s)^2 = -c^2 \overline{f'(p(s))}$. Taking the derivative with respect to $s$ yields

$$f''(p(s))p'(s)\overline{f(p(s))}^2 p'(s)^2 + f'(p(s))2\overline{f(p(s))}\,\overline{f'(p(s))p'(s)}p'(s)^2$$
$$+ f'(p(s))\overline{f(p(s))}^2 2p'(s)p''(s) = -c^2 \overline{f''(p(s))p'(s)}$$

which immediately yields

$$p''(s) = -\frac{c^2 \overline{f''(p(s))}\,\overline{p'(s)}^2}{2f'(p(s))\overline{f(p(s))}^2} - \frac{f''(p(s))p'(s)^2}{2f'(p(s))} - \frac{\overline{f'(p(s))}}{\overline{f(p(s))}}$$

after application of the identity $p'(s)\overline{p'(s)} = 1$. Higher derivatives may be obtained iteratively by differentiating with respect to $s$ and then making use of the identity $p'(s)\overline{p'(s)} = 1$. We then have the power series expansion

$$p(s) = \sum_{n=0}^{\infty} \frac{p^{(n)}(s_0)}{n!}(s-s_0)^n = p(s_0) + p'(s_0)(s-s_0)$$
$$-\frac{1}{2}\left(\frac{c^2\overline{f''(p(s_0))}\,\overline{p'(s_0)}^2}{2f'(p(s_0))\overline{f(p(s_0))}^2} + \frac{f''(p(s_0))p'(s_0)^2}{2f'(p(s_0))} + \frac{\overline{f'(p(s_0))}}{\overline{f(p(s_0))}}\right)(s-s_0)^2 + \cdots$$
$$= z_0 + \frac{if(z_0)\overline{f'(z_0)}}{|f'(z_0)||f(z_0)|}(s-s_0)$$
$$-\frac{1}{2}\left(\frac{-\overline{f''(z_0)}f'(z_0)}{2|f'(z_0)|^2} - \frac{f''(z_0)f(z_0)^2\,\overline{f'(z_0)}^3}{2|f'(z_0)|^4|f(z_0)|^2} + \frac{f(z_0)\overline{f'(z_0)}}{|f(z_0)|^2}\right)(s-s_0)^2 + \cdots$$

where

$$p'(s_0) = \frac{|f'(p(s_0))|ci}{f'(p(s_0))\overline{f(p(s_0))}} = \frac{|f'(z_0)|ci}{f'(z_0)\overline{f(z_0)}}.$$

Example: Let us consider the case $f(z) = az + b$, with $a \neq 0$ and $z_0 \neq -\dfrac{b}{a}$. The corresponding differential equation is

$$p'(s) = \frac{|a(ap(s)+b)|i}{a(ap(s)+b)}$$

Under the condition $p(s_0) = z_0$ this differential equation has the unique solution

$$p(s) = \left(z_0 + \frac{b}{a}\right)\exp\left(i\frac{1}{\left|z_0 + \dfrac{b}{a}\right|}(s - s_0)\right) - \frac{b}{a}.$$